\documentclass[a4paper,12pt]{article}

\usepackage[T1]{fontenc}
\usepackage[utf8]{inputenc}
\usepackage[english]{babel}
\usepackage{indentfirst}
\usepackage{amsfonts,amssymb,amsmath,amsthm}
\usepackage{mathtools}
\usepackage{array}
\usepackage{pdflscape}
\usepackage{float}
\usepackage{url}
\usepackage{hyperref}

\usepackage{geometry}
\newtheorem{theorem}{Theorem}
\newtheorem{proposition}{Proposition}
\newtheorem{definition}{Definition}

\title{The Borsuk Problem for Subsets of the Vertices\\ of the Boolean Cube in Dimensions 10 and 11}
\author{Igor Batmanov, Vsevolod Voronov}
\date{}

\begin{document}

\maketitle

\begin{abstract}
In the papers \cite{Ziegler99,Goldstein12} it was previously shown that any subset of the Boolean cube $S \subset \{0,1\}^n$ for $n \leq 9$ can be partitioned into $n+1$ parts of smaller diameter, i.e., the Borsuk conjecture holds for such subsets. In this paper, it is shown that this is also true for $n=10$ and $n=11$; however, the complexity of the computational verification increases significantly. In order to reduce the search tree, several heuristics were developed; the total computation time was less than two hours on a single core. The SAT solver \textbf{kissat} was used to cut off the search branches.
\end{abstract}

\section{Introduction}

This work is devoted to a particular case of the Borsuk problem on partitioning bounded subsets of Euclidean space $A \subset \mathbb{R}^n$ into parts of strictly smaller diameter. Let $k$ be the smallest number of parts for which such a partition exists,
\[
    A = A_1 \bigsqcup \dots \bigsqcup A_k, \quad \operatorname{diam} A_i < \operatorname{diam} A.
\]
Then $k$ is called the Borsuk number of the set $A$ and is denoted by $b(A)$. Furthermore, let
\[
    b(n) = \max b(A), \quad A \in \mathcal{A}_n,
\]
where $\mathcal{A}_n$ is the set of all bounded subsets of $\mathbb{R}^n$.

A similar problem can be formulated for subsets of a given (not necessarily bounded) set $M \subset \mathbb{R}^n$ and for an arbitrary metric space. In this paper, we consider the case of the Boolean cube $M = \{0,1\}^n$.

Counterexamples to the Borsuk conjecture in the form of subsets of $\{0,1\}^n$ are known for dimensions $n \geq 560$ \cite{Rai2,Rai3}. The smallest dimension for which a counterexample of a different type (a two-distance graph) has been constructed is currently $64$ \cite{Bond14,Jen14}.

The Borsuk conjecture for subsets of $\{0,1\}^n$ in low dimensions was previously studied in \cite{Ziegler99,Goldstein12}, where it was shown that such subsets can be partitioned into $n+1$ parts of smaller diameter for $n \leq 9$. In dimension $n=10$, the problem becomes qualitatively more complex. In this paper, we demonstrate that for $n=10$ and $n=11$ there is also no counterexample of the form $A \subset \{0,1\}^n$.

Unlike the approach in \cite{Goldstein12}, the search is not conducted for specific sets of vertices satisfying the diameter constraint but for a certain family of supersets. A universal covering system for a discrete set is considered, i.e., a family of subsets $\mathcal{T}$ such that any subset $V$ of a given diameter is mapped by some isometry into a subset of one of the sets in $\mathcal{T}$.

The main result of this work is the following.

\begin{theorem}
For any subset $V \subset \{0,1\}^n$ with $n \in \{10, 11\}$ we have $b(V) \leq n+1$, i.e., the Borsuk conjecture holds for subsets of the vertices of the Boolean cube in dimensions~10 and~11.
\end{theorem}

\section{Preliminary}

We introduce the notation for a distance graph with a given vertex set $V \subset \{0,1\}^n$,
\[
G(V; k) = (V,E), \qquad E = \{(u,v): \; \|u-v\|=k\}.
\]
Here and throughout, $\|v\| = \sum_i |v_i|$, and $\|u-v\|$ denotes the Hamming distance between $u$ and $v$. We define
\[
 G_{n,k} := G(\{0,1\}^n; k).
\]
For brevity, we denote
\[
    \chi\left(V; k\right) = \chi\left(G(V; k)\right)
\]
as the chromatic number of the graph $G(V; k)$.

Suppose that the Borsuk conjecture is false for some subset $V \subset \{0,1\}^n$ of diameter $k$. Then there is no proper coloring of the vertices of $G(V; k)$ using $k+1$ colors, and the following inequalities hold:
 \[
 \chi\left(G_{n, k}\right) \geq \chi\left(G(V; k)\right) \geq n+2.
 \]
Consequently, the first, relatively weak condition that allows establishing the validity of the conjecture for a given $n$ is the following.

\begin{proposition}
If $\forall k \; \chi(G_{n,k}) \leq n+1$, then the Borsuk conjecture holds for subsets of $\{0,1\}^n$.
\label{gnk_basic}
\end{proposition}

The study of the properties of the graphs $G_{n,k}$ is somewhat simplified by the following statements \cite{Ziegler99}.

\begin{proposition}
If $k$ is odd, then $\chi(G_{n,k})=2$.
\end{proposition}

\begin{proposition}
If $k$ is even, then $G_{n,k}$ consists of two isomorphic connected components, each of which contains only vertex pairs at even distances from each other.
\end{proposition}

It is easy to see that Propositions~1--3 imply the validity of the Borsuk conjecture for subsets of $\{0,1\}^n$ when $n \in \{4, 7, 8\}$, since in these cases there exist $(k+1)$-colorings of $G_{n,k}$ \cite{Ziegler99}. The cases $k=2$ for $5\leq n \leq 6$, $k \in \{2,4,6\}$ for $n=9$, and $k \in \{2,4,6,8\}$ for $n \in \{10, 11\}$ require a more detailed examination.

Let $S \subset V \subseteq \{0,1\}^n$ be a certain set of vertices. We define the operation of trimming with respect to the $k$-neighborhoods of the set $S$:
 \[
 \operatorname{Trim}(V, k, S)  = \{v \in V \, : \; \forall s \in S: \; \|v-s\| \leq k \}.
 \]

We denote by $B_{n,k}$ a certain subgraph of $G_{n,k}$ satisfying the conditions
 \[
 \forall u,v \in V(B_{n,k}) \quad \|u-v\| \leq k,
 \]
 \[
 \forall u \in \{0,1\}^n \setminus V(B_{n,k}) \quad \operatorname{diam}(V(B_{n,k}) \cup \{u\}) > k.
 \]
Let $\mathcal{B}_{n,k}$ be the set of all such subgraphs. Let us denote by $\mathcal{I}_n = \operatorname{Isom}(\{0,1\}^n)$ the symmetry group of the $n$-dimensional cube. Clearly, $|\mathcal{I}_n| = 2^n \cdot n!$. We are interested in considering subsets of the vertices of the cube, distinct up to symmetries.

\begin{definition}
We call a universal $(n,k)$-cover in $G_{n,k}$ (for brevity, simply an $(n,k)$-cover) a set $C_n \subset \{0,1\}^n$ such that
 \[
  \forall B \in \mathcal{B}_{n,k} \quad \exists g \in \mathcal{I}_n: \quad g \cdot B \subseteq C_n.
 \]
\end{definition}

Since for any $k$ the group $\mathcal{I}_n$ is transitive on the set of edges of $G_{n,k}$, we have the following.

\begin{proposition}
Let $u = (1, \dots , 1, 0, \dots, 0)$, with $\|u\|=k$. Then $T_{n,k,2} = \operatorname{Trim}(\{0,1\}^n, k, \{0_n, u\})$ is an $(n,k)$-cover.
\end{proposition}

This statement allows us to prove the $(0,1)$-Borsuk conjecture for $n = 9$ using computer calculations.

\begin{proposition}
$\chi\left(T_{9,k,2}; k\right) \leq 10, \; 1\leq k\leq 9$.
\end{proposition}

\begin{proposition}
$\chi\left(T_{10,k,2}; k\right) \leq 11$ for $k \in \{2, 8\}$.
\end{proposition}

\begin{proposition}
$\chi\left(T_{11,k,2}; k\right) \leq 12$ for $k \in \{2, 8\}$.
\end{proposition}

Propositions~5--7 were verified by computer. However, for $n \in \{10, 11\}$ and $k \in \{4, 6\}$ the cover $T_{n,k,2}$ turns out to be too large for the direct computation of the chromatic number. Instead, a family of covering sets is introduced.

\section{Covers and covering systems}

To prove the Borsuk conjecture for fixed $n$ and $k$, one must show that \emph{every} subset $V \subset \{0,1\}^n$ of diameter $k$ admits a proper coloring in $n+1$ colors. The key idea is to cover the space of all such subsets by a finite family of supersets, for each of which colorability can be verified by computer.

By analogy with the constructions used in the estimates for the Borsuk problem in the case of compact subsets of $\mathbb{R}^n$, we introduce the following.

\begin{definition}
A universal covering $(n,k)$-system is a family of sets $\mathcal{C}_{n,k}$ such that
 \[
 \forall C \in \mathcal{C}_{n,k} \quad C \subset \{0,1\}^n,
 \]
 \[
 \forall B \in \mathcal{B}_{n,k} \quad \exists g \in \mathcal{I}_n: \quad \exists C \in \mathcal{C}_{n,k}: \quad g \cdot B \subseteq C.
 \]
\end{definition}

It is obvious that if
\[
      \chi(C; k) \leq n+1, \quad \forall C \in \mathcal{C}_{n,k},
\]
then the Borsuk conjecture is proven in this case.

\subsection*{Case $n=10$, $k=4$}

Let us illustrate the construction of a covering system with a relatively simple example. For brevity, we write the coordinates of a vector as a string without using commas. Let
\[
u_1 = (0000001111),\; u_2 = (0000110011),\; v_1 = (0000010111),
\]
\[
v_2 = (0000100111),\; \ldots\,,\; v_6 = (1000000111),\; w = (0000011011);
\]
 \[
 U_1 = \{0_{10},u_1,u_2\}, \quad
 U_2 = \{0_{10},u_1,v_1,\dots, v_6\}, \quad
 U_3 = \{0_{10},u_1,v_1,w\},
 \]
 \[
 W = \operatorname{Trim}(\{0,1\}^{10},2,\{0_{10}\}).
 \]

\begin{proposition}
The family
 \[
 \mathcal{C}^*_{10,4} = \{ \operatorname{Trim}(\{0,1\}^{10},4, U_1), \; U_2 \cup W, \; U_3 \cup W \}
 \]
is a covering $(10,4)$-system.
\end{proposition}

\noindent\emph{Proof.} If the distance graph of a set $B$ contains a triangle $(4,4,4)$, then $B$ is covered by the first element of the system. Otherwise, the number of vertices adjacent to $0$ is limited, and all other vertices are at distance at most $2$ from $0$. If the vertices with four unit coordinates have zero instead of one in the same position compared to $u_1$, we obtain the set $U_2$; if in different ones, we get $U_3$. \hfill$\square$

\begin{proposition}
For every $C \in \mathcal{C}^*_{10,4}$ we have $\chi(C; 4) \leq 11$.
\end{proposition}

\noindent\emph{Proof.} Computer verification. Analogous statements hold for all pairs $(n,k)$ considered in this work.

\begin{proposition}
For each pair
 \[
 (n, k) \in \{(10,4),\; (10,6),\; (11,4),\; (11,6)\}
 \]
there exists a covering $(n,k)$-system $\mathcal{C}^*_{n,k}$ such that $\chi(C; k) \leq n+1$ for all $C \in \mathcal{C}^*_{n,k}$.
\end{proposition}

\noindent\emph{Proof.} Computer verification (see the section ``Results'').

\section{Classification of complete subgraphs in $G_{n,6}$}

The case $k=6$ is considerably more complicated than those considered above. The direct computation of the chromatic number of the graph $G_{n,6}$ for $n \geq 10$ is impossible because of the size of the graph. The construction of a covering system for $k=6$ is based on the classification of complete subgraphs (cliques) in $G_{n,6}$ up to equivalence.

\begin{definition}
Two sets $S_1, S_2 \subset \{0,1\}^n$ are called \emph{equivalent} if there exists an isometry $g \in \mathcal{I}_n$ such that $g \cdot S_1 = S_2$.
\end{definition}

If the presence of a subgraph $K_m$ is assumed, one can fix its realization and discard all vertices located at distance greater than $k$ from at least one of the vertices of $K_m$. Moreover, if the absence of $K_{m+1}$ is assumed, then the vertices that would form a $K_{m+1}$ together with the already chosen ones can be forbidden in order to simplify the computations.

Since the group $\mathcal{I}_n$ is transitive on the vertices, any clique is equivalent to a clique containing $0_n$. Therefore it is sufficient to classify cliques containing $0_n$. The stabilizer of $0_n$ in $\mathcal{I}_n$ coincides with the group of coordinate permutations $S_n$, which motivates the following definition.

\begin{definition}
Let $S = \{v_0 = 0_n, v_1, \ldots, v_{m-1}\} \subset \{0,1\}^n$ be a set of vertices containing $0_n$. The \emph{coordinate profile} of the set $S$ is the vector $\sigma(S) = (s_1^*, \ldots, s_n^*)$, where $(s_1^*, \ldots, s_n^*)$ is the permutation of the vector $(s_1, \ldots, s_n)$ in non-decreasing order, and
\[
s_j = \sum_{v \in S} v_j
\]
is the number of vertices of the set $S$ having a one in the $j$-th coordinate. A Young diagram with $n$ columns of heights $s_1^*, \ldots, s_n^*$ is naturally associated with the coordinate profile.
\end{definition}

Since a coordinate permutation only permutes the components of the vector $(s_1, \ldots, s_n)$, the coordinate profile $\sigma(S)$ is an invariant of the equivalence class (among sets containing $0_n$).

\textbf{Remark.} The coordinate profile is a necessary but not a sufficient condition for equivalence. However, for cliques $K_m$ in $G_{n,6}$ with $m \leq 5$ and $n \leq 12$, the profile completely determines the equivalence class.

Any two instances of $K_2$ and $K_3$ in $G_{n,6}$ are equivalent; however, for $K_4$ and $K_5$ there exist several equivalence classes, distinguished by the coordinate profile.

\subsection*{Equivalence classes of $K_4$}

There are $4$ equivalence classes of $K_4$ in $G_{n,6}$:
\begin{itemize}
 \item For $n \geq 9$: $K^{(1)}_4 \sim (0_9,\;
 \texttt{111\,111\,000},\;
 \texttt{111\,000\,111},\;
 \texttt{000\,111\,111})$,\quad $\sigma = (2^9)$.

 \item For $n \geq 10$: $K^{(2)}_4 \sim (0_{10},\;
 \texttt{111\,111\,000\,0},\;
 \texttt{111\,000\,111\,0},\;
 \texttt{010\,101\,101\,1})$,\quad $\sigma = (1^3, 2^6, 3)$.

 \item For $n \geq 11$: $K^{(3)}_4 \sim (0_{11},\;
 \texttt{111\,111\,000\,00},\;
 \texttt{111\,000\,111\,00},\;
 \texttt{101\,010\,010\,11})$,\quad $\sigma = (1^6, 2^3, 3^2)$.

 \item For $n \geq 12$: $K^{(4)}_4 \sim (0_{12},\;
 \texttt{111\,111\,000\,000},\;
 \texttt{111\,000\,111\,000},\;
 \texttt{111\,000\,000\,111})$,\quad $\sigma = (1^9, 3^3)$.
\end{itemize}

\subsection*{Equivalence classes of $K_5$}

There are $5$ equivalence classes of $K_5$ in $G_{n,6}$:
\begin{itemize}
 \item For $n \geq 10$: $K^{(1)}_5 \sim (0_{10},\;
 \texttt{111\,111\,000\,0},\;
 \texttt{111\,000\,111\,0},\;
 \texttt{010\,101\,101\,1},\;
 \texttt{100\,011\,011\,1})$,\quad $\sigma = (2^6, 3^4)$.

 \item For $n \geq 11$: $K^{(2)}_5 \sim (0_{11},\;
 \texttt{111\,111\,000\,00},\;
 \texttt{111\,000\,111\,00},\;
 \texttt{010\,101\,101\,10},\;
 \texttt{010\,110\,110\,01})$,\quad $\sigma = (1^2, 2^6, 3^2, 4)$.

 \item For $n \geq 11$: $K^{(3)}_5 \sim (0_{11},\;
 \texttt{111\,111\,000\,00},\;
 \texttt{111\,000\,111\,00},\;
 \texttt{010\,101\,101\,10},\;
 \texttt{100\,101\,110\,01})$,\quad $\sigma = (1^3, 2^3, 3^5)$.

 \item For $n \geq 12$: $K^{(4)}_5 \sim (0_{12},\;
 \texttt{111\,111\,000\,000},\;
 \texttt{111\,000\,111\,000},\;
 \texttt{101\,010\,010\,101},$\\
 $\texttt{101\,100\,100\,110})$,\quad $\sigma = (1^4, 2^6, 4^2)$.

 \item For $n \geq 12$: $K^{(5)}_5 \sim (0_{12},\;
 \texttt{111\,111\,000\,000},\;
 \texttt{111\,000\,111\,000},\;
 \texttt{101\,010\,010\,101},$\\
 $\texttt{110\,010\,100\,110})$,\quad $\sigma = (1^5, 2^3, 3^3, 4)$.
\end{itemize}

\medskip

The classification described above determines a decomposition of the problem into subcases. In each subproblem a separate covering system is constructed by the algorithm described below.

For $n=10$ there are $5$ subproblems, corresponding to the subgraphs $K_2$, $K_3$, $K^{(1)}_4$, $K^{(2)}_4$, $K_5$. For $n=11$ the types $K^{(3)}_4$, $K^{(1)}_5$, $K^{(2)}_5$, $K^{(3)}_5$ are added, which gives $8$ subproblems. The types $K^{(4)}_4$, $K^{(4)}_5$, $K^{(5)}_5$ arise only for $n \geq 12$ and are given for the completeness of the classification.

\section{Construction of the covering system for $k=6$}

We introduce an extension of the operation $\operatorname{Trim}$ that takes forbidden constructions into account:
\begin{multline*}
\operatorname{Trim2}(V, k, S, F) = \{v \in V : \forall s \in S\; \|v-s\| \leq k \\
\;\wedge\; \nexists f \in F: f \subset \{v\} \cup S \}.
\end{multline*}
Here $S$ is the set of mandatory vertices, and $F$ is the set of forbidden constructions (complete subgraphs whose presence is excluded in this subproblem).

Consider a set of mandatory vertices $U$, a set of forbidden vertices $C$, and a set of forbidden constructions $F$. If it is possible to verify that
\[
\chi(\operatorname{Trim2}(\{0,1\}^n \setminus C, k, U, F); k) \leq n+1,
\]
we call the configuration $(U, C, F)$ \emph{successful}.

The algorithm constructs a set of successful configurations $\mathcal{U}^*$ by recursive enumeration. Let a set of vertices $U_t = \{u_1, \dots, u_t\}$ and a set of forbidden vertices $C_g = \{c_1, \dots, c_g\}$ be selected. If the configuration is not successful, a vertex
\[u_{t+1} \in \operatorname{Trim2}(\{0,1\}^n \setminus C_g, k, U_t, F)\]
is selected, and two branching options are considered:
\begin{enumerate}
    \item the vertex $u_{t+1}$ is added to the mandatory ones: $U_{t+1} = U_t \cup \{u_{t+1}\}$;
    \item the vertex $u_{t+1}$ is forbidden: $C_{g+1} = C_g \cup \{u_{t+1}\}$.
\end{enumerate}

If a set of successful configurations $\mathcal{U}^*$ is found such that any set from $\mathcal{B}_{n,k}$ contains at least one of them, then
  \[
  \mathcal{T}_{n,k} = \{\operatorname{Trim2}(\{0,1\}^n \setminus C_U, k, U, F_U) \mid U \in \mathcal{U}^* \}
 \]
is a covering $(n,k)$-system, and for each element the chromatic number condition is satisfied. Each leaf of the search tree in which a coloring is successfully found generates one element of the covering system.

The set of vertices for the enumeration is split into subsets $M_1, M_2, \ldots, M_f$; the enumeration proceeds sequentially. To check colorability, the SAT solver \textbf{kissat} is used with a time limit (10--120\,s). If a coloring is not found within the allotted time, the configuration is considered not \emph{successful}.

\section{Pseudocode}

Each element of the covering system has the form $\operatorname{Trim2}(\{0,1\}^n \setminus C_U, k, U, \mathcal{F}_U)$. In terms of the algorithm:
\begin{itemize}
    \item \textbf{config} --- the set of mandatory vertices $U$;
    \item \textbf{rest} $\cup$ \textbf{setToCheck} --- the remaining vertices from $\operatorname{Trim2}(\{0,1\}^n \setminus C_U, k, U, \mathcal{F}_U)$;
    \item \textbf{a successful configuration} (the graph is colored in $n+1$ colors) --- an element of the covering system $C \in \mathcal{T}_{n,k}$;
    \item each leaf of the search tree in which a coloring is found generates an element of the covering system.
\end{itemize}

\textbf{Global variables}
\begin{itemize}
    \item $WeakSymmetries$ --- the set of weak symmetries preserving a given construction $X$. The transformations include coordinate permutations and mirror reflections (XOR operations). The name ``weak'' reflects the fact that this is a subset of all symmetries, used to speed up the verification.
    \item $ProcessedHashes$ --- the set of hashes of already processed configurations (initially empty).
    \item $NotColored$ --- the set of configurations for which a coloring could not be found within the allotted time (initially empty).
\end{itemize}

\hrulefill

\textbf{Main procedure}

\begin{flushleft}
\textbf{Require:} $X_{assumptions}$, $n$, $\mathcal{F}_{forbidden}$, $M_1, M_2, \dots$\\[0.4ex]
$WeakSymmetries \leftarrow \mathrm{calcWeakSymmetries}(X_{assumptions}, n)$\\
$ProcessedHashes \leftarrow \emptyset$;\quad $NotColored \leftarrow \emptyset$\\
$rest \leftarrow \mathrm{Trim2}(\{0,1\}^n, k, X_{assumptions}, \mathcal{F}_{forbidden})$\\
$Branches(X_{assumptions}, M_1, rest)$\\[0.4ex]
\textbf{if} $NotColored \neq \emptyset$ \textbf{then}\\
\hspace*{1.5em}\textbf{for} $(config, restVertices) \in NotColored$ \textbf{do}\\
\hspace*{3em}$Branches(config, M_2, restVertices)$\\
\hspace*{1.5em}\textbf{end for}\\
\textbf{end if}\\[0.4ex]
\textbf{if} $NotColored \neq \emptyset$ \textbf{then return} ``Failed to color all configurations''\\
\textbf{else return} ``All configurations are colored in $n+1$ colors''\\
\textbf{end if}
\end{flushleft}

\hrulefill

\textbf{wasBefore(config)}

\begin{flushleft}
\textbf{Require:} $config$ --- the current configuration of vertices\\[0.4ex]
$h \leftarrow \mathrm{hash}(config)$\\
\textbf{if} $h \in ProcessedHashes$ \textbf{then return} $true$\quad\textbf{end if}\\
\textbf{for} $T \in WeakSymmetries$ \textbf{do}\\
\hspace*{1.5em}$h_T \leftarrow \mathrm{hash}(T(config))$\\
\hspace*{1.5em}$ProcessedHashes \leftarrow ProcessedHashes \cup \{h_T\}$\\
\textbf{end for}\\
\textbf{return} $false$
\end{flushleft}

\hrulefill

\textbf{Branches(config, setToCheck, rest)}

\begin{flushleft}
\textbf{if} $wasBefore(config)$ \textbf{then return}\quad\textbf{end if}\\
\textbf{if} $isColorable(config, setToCheck \cup rest)$ \textbf{then return}\quad\textbf{end if}\\
\textbf{if} $\neg\, isColorable(config, rest)$ \textbf{then}\\
\hspace*{1.5em}$NotColored \leftarrow NotColored \cup \{(config, rest)\}$;\quad\textbf{return}\\
\textbf{end if}\\
\textbf{for} $newVertex \in setToCheck$ \textbf{do}\\
\hspace*{1.5em}\textbf{if} $canAddRestricted(config, newVertex)$ \textbf{then}\\
\hspace*{3em}$Branches(config \cup \{newVertex\},\ setToCheck[newVertex\!+\!1\!:],\ rest)$\\
\hspace*{1.5em}\textbf{end if}\\
\textbf{end for}
\end{flushleft}

\hrulefill

\section{Key optimizations}

\begin{itemize}
    \item \textbf{wasBefore}: for each configuration a hash is computed; then the weak symmetries ($WeakSymmetries$) --- a subset of isometries preserving the starting construction --- are applied. If the hash of a transformed configuration has already been encountered, the branch is cut off.
    \item \textbf{canAddRestricted}: when a vertex is added, it is checked whether it forms a forbidden construction together with the already selected vertices.
    \item The size of the starting set $M_1$ is chosen so that the set of not-colored configurations $NotColored$ remains empty.
\end{itemize}

If the algorithm terminates successfully ($NotColored = \emptyset$), the constructed set $\mathcal{U}^*$ generates a covering $(n,k)$-system, which proves the conjecture.

\section{Results}

The algorithm described above was applied to all nontrivial pairs $(n, k)$: for $n=10$ to the cases $k \in \{2,4,6,8\}$, and for $n=11$ to the cases $k \in \{2,4,6,8\}$. For $k=6$ the enumeration was carried out over all subproblems determined by the clique classification: $5$ subproblems for $n=10$ and $8$ subproblems for $n=11$ (see Tables~\ref{tab:n10} and~\ref{tab:n11}).

In none of the runs did an uncolored configuration arise ($NotColored = \emptyset$). Thereby, for each pair $(n, k)$ a covering system was constructed, all elements of which admit a coloring in $n+1$ colors, which proves the validity of the Borsuk conjecture for subsets of $\{0,1\}^n$ for $n=10$ and $n=11$.

Without the use of parallel computations, on a single core with a frequency of $3$\,GHz, the total computation time was less than two hours. The source code is available in the repository \cite{github}.

In Tables~\ref{tab:n10} and~\ref{tab:n11}, the values in the columns TLeaves, OLeaves, Skipped are given in the format ``depth:count''. TLeaves denotes the terminal leaves of the search tree in which a coloring was successfully found. OLeaves denotes the leaves in which a coloring was not found within the allotted time (configurations from $NotColored$, processed a second time with the set $M_{i+1}$ and with the remaining vertices of $M_i$ additionally forbidden). Skipped denotes the configurations skipped thanks to the symmetry cache ($wasBefore = true$).

\begin{landscape}
\begin{table}[H]
    \centering
    \begin{tabular}{|c|c|c|c|c|c|c|c|}
    \hline
        No & S & Fbd & $M_1$ & TLeaves & OLeaves & Skipped & Time  \\
    \hline
      N10k2 & $K_2$ & no & 100 & 2:1 & 0 & 0 &  0min10s  \\
    \hline
      N10k4 & $K_2$ & no & 100 & 3:2 & 2:1 & 3:98 &  21s  \\
    \hline
      N10k8 & $K_2$ & no & 100 & 2:1 & 0 & 0 &  0min10s  \\
    \hline
      N10k6 & $K_2$ & $K_3$ & 150 & 4:10 & 2:1, 3:1 & 3:149, 4:97 & 33s  \\
    \hline
      N10k6 & $K_3$ & $K_4$ & 100 & 5:13, 6:72 & 3:1, 4:2, & 4:98, 5:136, &  2min1s  \\
       & & & &  & 5:4, 6:1 & 6:108 &   \\
    \hline
      N10k6 & $K^{(1)}_4$ & $K_5$ & 100 & 6:2, 7:50 & 4:1, 5:1,  & 5:99, 6:81,  &  1min55s  \\
      & & & & & 6:5, 7:1& 7:255& \\
    \hline
      N10k6 & $K^{(2)}_4$ & $K_5$ & 100 & 6:83, 7:636 & 4:1, 5:5,& 5:95, 6:202, &  7min19s  \\
      & & & & &  6:19, 7:8&  7:440& \\
    \hline
      N10k6 & $K_5$ & no & 100 & 7:42, 8:665,  & 5:1, 6:5, & 6:95, 7:355,  &  32min  \\
      & & & & 9:2701,&  7:20, 8:63,  & 8:794, 9:1235, & 53s\\
      & & & &  10:2088 & 9:72, 10:28& 10:527  & \\
    \hline
    \end{tabular}
    \caption{Case enumeration for $n=10$ at different $k$}
    \label{tab:n10}
\end{table}
\end{landscape}

\begin{landscape}
\begin{table}[H]
    \centering
    \begin{tabular}{|c|c|c|c|c|c|c|c|}
    \hline
        No & S & Fbd & $M_1$ & TLeaves & OLeaves & Skipped & Time  \\
    \hline
      N11k2 & $K_2$ & no & 100 & 2:1 & 0 & 0 &  1min58s  \\
    \hline
      N11k4 & $K_2$ & no & 100 & 3:2 & 2:1 & 3:98 &  1min50s  \\
    \hline
      N11k8 & $K_2$ & no & 100 & 2:1 & 0 & 0 &  1min55s  \\
    \hline
      N11k6 & $K_2$ & $K_3$ & 200 & 4:10 & 2:1, 3:1 & 3:199, 4:114 & 2m22s  \\
    \hline
      N11k6 & $K_3$ & $K_4$ & 201 & 4:1, 5:24, & 3:1, 4:3, & 4:197, 5:353&  6min  \\
      & & & & 6:216& 5:8, 6:3& 6:298&15s \\
    \hline
      N11k6 & $K^{(1)}_4$ & $K_5,$ & 201 & 4:1 & 0 & 0 &  4min28s  \\
      & & $K^{(3)}_4,$& & & & & \\
      & & $K^{(2)}_4$& & & & & \\
    \hline
      N11k6 & $K^{(2)}_4$ & $K_5,$ & 201 & 5:1, 6:195, & 4:1, 5:7,  & 5:193, 6:643,  &  9min51s  \\
& &  $K^{(3)}_4$& & 7:446& 6:9, 7:5& 7:347& \\
    \hline
      N11k6 & $K^{(3)}_4$ & $K_5$ & 201 & 5:3, 6:62, & 4:1, 5:2, & 5:196, 6:233,&  8min35s  \\
    & & & & 7:635 &  6:11, 7:5&  7:517 & \\
    \hline
      N11k6 & $K^{(1)}_5$ & no & 216 & 6:4, 7:238, & 5:1, 6:7,  & 6:205, 7:840, &  16min41s  \\
    & & & & 8:907, 9:668 & 7:16, 8:12, 9:4&  8:1261, 9:267& \\
    \hline
      N11k6 & $K^{(2)}_5$ & no & 216 & 6:12, 7:843, & 5:1, 6:13,& 6:191, 7:964 &  11min18s  \\
    & & & & 8:126 &  7:6, 8:1 & & \\
    \hline
      N11k6 & $K^{(3)}_5$ & no & 216 & 6:13, 7:874 & 5:1, 6:13, 7:5 & 6:190, 7:966 &  11min24s  \\
    \hline
    \end{tabular}
    \caption{Case enumeration for $n=11$ at different $k$}
    \label{tab:n11}
\end{table}
\end{landscape}

\section{Conclusion}

The calculations presented above suggest that the existence of a distance graph with a chromatic number of $n+2$ on the vertices of the cube of dimension $n$ seems quite plausible even in relatively low dimensions. Such a counterexample may be found as a result of solving the independent set search problem in a prohibition graph, where edges are drawn between vertices at distances greater than $k$. However, the direct application of algorithms for the maximum independent set does not yield the desired result; the chromatic number of the constructed diameter graph turns out to be much smaller than the maximum possible (for example, the chromatic number does not exceed $7$ in dimension $16$). Apparently, it is at least necessary to add weights to the problem formulation, or, more likely, to develop a special modification of the algorithm.

Moreover, it is evident that a multitude of potential (or even actual) counterexamples can be constructed without any reasonable justification that the chromatic number of the diameter graph is indeed greater than or equal to $n+2$. Existing counterexamples are constructed precisely in such a way that the proof is feasible. Therefore, in order to have a chance of obtaining the required estimate, the construction must be redundant: if the diameter graph is difficult to color even with $n+k$ colors, $k \geq 2$, there is a chance that the computer verification of the absence of an $(n+1)$-coloring will finish in a reasonable time.


\begin{thebibliography}{99}

\bibitem{Kahn93}
Kahn J., Kalai G. A counterexample to Borsuk's conjecture // Bulletin of the American Mathematical Society. --- 1993. --- Vol.~29. --- No.~1. --- P.~60--62.

\bibitem{Nilli}
Nilli, A. (1994). On Borsuk's problem. Contemporary Mathematics, 178, 209--209.

\bibitem{Rai1}
Raigorodskii, A. M. (2006). The Borsuk problem. Moscow: MCCME, 130.

\bibitem{Rai2}
Raigorodskii, A. M. (2001). The Borsuk problem and the chromatic numbers of some metric spaces. Russian Mathematical Surveys, 56(1 (337)), 107--146.

\bibitem{Rai3}
Raigorodskii, A. M. (2007). Around the Borsuk conjecture. Contemporary Mathematics. Fundamental Directions, 23(0), 147--164.

\bibitem{Bond14}
Bondarenko, A. (2014). On Borsuk's conjecture for two-distance sets. Discrete and Computational Geometry, 51(3), 509--515.

\bibitem{Jen14}
Jenrich, T., Brouwer, A. E. (2014). A 64-dimensional counterexample to Borsuk's conjecture. The Electronic Journal of Combinatorics, P4-29.

\bibitem{Jen23}
Jenrich, T. (2023). Sub-25-dimensional counterexamples to Borsuk's conjecture in the Leech lattice? arXiv preprint arXiv:2305.06283.

\bibitem{Hin03}
Hinrichs, A., Richter, C. (2003). New sets with large Borsuk numbers. Discrete Mathematics, 270(1--3), 137--147.

\bibitem{Ziegler99}
Ziegler, G. M. (2001). Coloring Hamming graphs, optimal binary codes, and the 0/1-Borsuk problem in low dimensions. In Computational Discrete Mathematics (pp. 159--171). Springer, Berlin, Heidelberg.

\bibitem{Goldstein12}
Goldstein, V. B. (2012). On the Borsuk problem for (0,1)- and (-1,0,1)-polytopes in low-dimensional spaces. Transactions of the Moscow Institute of Physics and Technology, 4(1--13).

\bibitem{Kokkala18}
Kokkala, J., Östergård, P. (2018). The chromatic number of the square of the 8-cube. Mathematics of Computation, 87(313), 2551--2561.

\bibitem{Matthews17}
Matthews, G. L. (2017). Distance colorings of hypercubes from $\mathbb{Z}_2\mathbb{Z}_4$-linear codes. Discrete Applied Mathematics, 217, 356--361.

\bibitem{Saras20}
Saraswathi, S., Poobalaranjani, M. (2020). Exact 2-distance b-coloring of some classes of graphs. Malaya Journal of Matematik (MJM), 8(1), 195--200.

\bibitem{kissat}
Fleury, M., Biere, A. (2023). Mining definitions in Kissat with Kittens. Formal Methods in System Design, 1--24.

\bibitem{handbook}
Biere, A., Heule, M., van Maaren, H. (Eds.). (2009). Handbook of Satisfiability (Vol.~185). IOS Press.

\bibitem{github}
Batmanov, I., Voronov, V. BorsukResearch --- source code repository. \url{https://github.com/igorbat/BorsukResearch}

\end{thebibliography}
\end{document}